\documentclass[conference]{IEEEtran}

\IEEEoverridecommandlockouts
\usepackage{cite}

\ifCLASSINFOpdf
   \usepackage[pdftex]{graphicx}
   \DeclareGraphicsExtensions{.pdf,.jpeg,.png}
\else
  \usepackage[dvips]{graphicx}
  \DeclareGraphicsExtensions{.eps}
\fi

\usepackage{textcomp}
\usepackage{xcolor}

\usepackage{amsthm}
\usepackage{lipsum}
\usepackage{color}
\usepackage{algorithm}
\usepackage{algorithmic}
\usepackage{enumerate}

\usepackage{etoolbox}  
\makeatletter
\patchcmd{\algorithmic}{\addtolength{\ALC@tlm}{\leftmargin} }{\addtolength{\ALC@tlm}{\leftmargin}}{}{}
\usepackage{amssymb}

\ifCLASSINFOpdf
   \usepackage[pdftex]{graphicx}
\else
\fi
\usepackage{cite}
\usepackage[cmex10]{amsmath}
\usepackage{algorithmic}
\usepackage{array}
\usepackage{url}
\usepackage{makeidx}
\usepackage{verbatim}
\usepackage{subfigure}

\usepackage{multirow}
\newcommand{\beqn}{\begin{eqnarray}}
\newcommand{\eeqn}{\end{eqnarray}}
\newcommand{\bse}{\begin{subequations}}
\newcommand{\ese}{\end{subequations}}

\begin{document}

\title{Multi-UAV trajectory planning problem using  the difference of convex function programming

\author{Anh Phuong Ngo$^1$, Christian Thomas$^2$, Ali Karimoddini$^1$ and Hieu T. Nguyen$^1$\\ 
$^1$\textit{Dept. of Electrical \& Computer Eng., North Carolina A\&T State University, Greensboro, NC 27411, USA}\\
$^2$\textit{Dept. of Flight Test, Lockheed Martin Corporation, Fort Worth, TX 76108, USA}\\
\textit{ango1@aggies.ncat.edu, christian.thomas@lmco.com, \{akarimod, htnguyen1\}@ncat.edu}

\thanks{
This work was supported by DOE Sandia National Laboratories under contract 281247. This paper has been accepted for presentation at the 62nd IEEE Conference on Decision and Control (CDC 2023).
}
}%
}

\maketitle
\begin{abstract}
The trajectory planning problem for a swarm of multiple UAVs is known as a challenging nonconvex optimization problem, particularly due to a large number of collision avoidance constraints required for individual pairs of UAVs in the swarm.
In this paper, we tackle this nonconvexity by leveraging the difference of convex function (DC) programming.
We introduce the slack variables to relax and reformulate the collision avoidance conditions and employ the penalty function term to equivalently convert the problem into a DC form. 
Consequently, we construct a penalty DC algorithm in which we sequentially solve a set of convex optimization problems obtained by linearizing the collision avoidance constraint.
The algorithm iteratively tightens the safety condition and reduces the objective cost  of the planning problem and the additional penalty term.
Numerical results demonstrate the effectiveness of the proposed approach in planning a large number of UAVs in congested space.

\end{abstract}



{ 
\footnotesize
\section*{Nomenclature}
\addcontentsline{toc}{section}{Nomenclature}
\subsection{Set and Indices}
\begin{IEEEdescription}[\IEEEusemathlabelsep\IEEEsetlabelwidth{$V_1,V_2,V_3$}]
\raggedright 
\item[$\mathcal{N}, i$] Set and index of vehicles, $i \in \mathcal{N}$
\item[$\mathcal{T}, k$] Set and index of time steps, $k=1,2, \dots, K \in \mathcal{T}$
\item[$\mathcal{S}$] Set of initial states, including starting position $\mathcal{S}\big(x_i^s,y_i^s,z_i^s \big)$, starting velocity $\mathcal{S}\big(v^x_i,v^y_i,v^z_i \big)$, and starting force $\mathcal{S}\big(f{^x_i,f^y_i,f^z_i} \big)$ of vehicle $i$
\item[$\mathcal{G}$] Set of goal states, including goal position $\mathcal{G}\big(x_i^g,y_i^g,z_i^g \big)$, goal velocity $\mathcal{G}\big(v^x_i,v^y_i,v^z_i \big)$, and goal force $\mathcal{G}\big(f{^x_i,f^y_i,f^z_i} \big)$ of vehicle $i$
\end{IEEEdescription}

\subsection{Parameters}
\begin{IEEEdescription}[\IEEEusemathlabelsep\IEEEsetlabelwidth{$V_1,V_2,V_3$}]
\raggedright 
\item[$\overline{x}, \underline{x}$] Upper/lower limits of $x$-coordinate that vehicles can reach
\item[$\overline{y}, \underline{y}$] Upper/lower limits of $y$-coordinate that vehicles can reach
\item[$\overline{z}, \underline{z}$] Upper/lower limits of $z$-coordinate that vehicles can reach
\item[$\overline{V_i}, \underline{V_i}$] Upper/lower limits of velocity of vehicle $i$
\item[$\overline{F_i}, \underline{F_i}$] Upper/lower limits of force of vehicle $i$
\item[$d$] Minimum distance among two vehicles to avoid a collision
\item[$\mathbf{A}_{i}$] State-space matrix of vehicle $i$
\item[$\mathbf{B}_{i}$] Input matrix of vehicle $i$
\item[$\rho^f$] Penalty for $O_f$ of objective function
\item[$\rho^k$] Penalty for $O_g$ of objective function
\item[$\tau, \mu, \epsilon$] Parameters used in DCA
\end{IEEEdescription}

\subsection{Variables}
\begin{IEEEdescription}[\IEEEusemathlabelsep\IEEEsetlabelwidth{$V_1,V_2,V_3$}]
\raggedright
\item[$x_{i,k}$] $x$-coordinate of position of vehicle $i$ at time step $k$
\item[$y_{i,k}$] $y$-coordinate of position of vehicle $i$ at time step $k$
\item[$z_{i,k}$] $z$-coordinate of position of vehicle $i$ at time step $k$
\item[$v^x_{i,k}$] $x$-component of velocity vector of vehicle $i$ at time step $k$
\item[$v^y_{i,k}$] $y$-component of velocity vector of vehicle $i$ at time step $k$
\item[$v^z_{i,k}$] $z$-component of velocity vector of vehicle $i$ at time step $k$
\item[$\Vec{V}_{i,k}$] velocity vector of vehicle $i$ at time step $k$
\item[$f^x_{i,k}$] $x$-component of force vector of vehicle $i$ at time step $k$
\item[$f^y_{i,k}$] $y$-component of force vector of vehicle $i$ at time step $k$
\item[$f^z_{i,k}$] $z$-component of force vector of vehicle $i$ at time step $k$
\item[$\Vec{F}_{i,k}$] force vector of vehicle $i$ at time step $k$
\end{IEEEdescription}
}

\section{Introduction}

The trajectory planning problem aims at finding an optimal solution of the trajectory for a single aircraft or a group of aircrafts to travel from a given starting state over a map of the environment to a goal state. 
Mixed-integer linear programming (MILP) is the standard method used to solve the trajectory generation problem for many decades \cite{MIP_review}.
MILP is a powerful optimization method that allows inclusion of integer variables and discrete logic of linearizaion for non-convex constraints in a continuous linear trajectory optimization \cite{UAVcomplexityMIT2002,taskallocation2001,2008robust}. 
These mixed-integer and continuous variables can be used to model logical constraints such as obstacle avoidance and vehicle separation, while the dynamic and kinematic settings of the aircrafts are bounded in continuous constraints.
Concurrently, the magnitudes of velocity and force vectors are modeled by the spherical geometry-based sampling approximation technique for a 3-D environment, or the edges of an N-sided polygon approximation technique for a 2-D environment \cite{2008robust}.
To this extent, the MILP method uses many auxiliary variables and constraints to formulate the trajectory optimization problem.

Recent improvements in aircraft's capabilities, especially for unnamed aerial vehicles (UAVs), facilitate them to carry out longer and more complex missions in dynamic environments.
Moreover, as more vehicles and more targets are involved in a mission, the size of the trajectory optimization problem based on MILP increases exponentially. 
Consequently, the computation time of the problem to obtain the optimal solution becomes much more expensive.
Convex optimization methods can handle well the conic constraints such as bounds on the magnitude of velocity and force vectors without incorporating the approximation techniques \cite{Convex_possibility_caltech, UAVcomplexityMIT2002}. 

The most key challenge in solving trajectory optimization models with convex cost functions and affine vehicle dynamics is that we often encounter the nonconvex collision avoidance requirement \cite{multi_air_weather_AIAA, realtimeUAV_Dynamic_Chin}.
This nonconvex requirement is enforced for all individual pairs of UAVs in the swarm, thus making the problem computationally challenging. 
This research proposes the use of the difference of convex function (DC) programming \cite{DCA_tao96} to tackle the nonconvexity of the planning problem for a swarm of a large number of UAVs.
First, we relax the collision avoidance constraints by slack variables and add the sum of slack variables as a penalty function to the original objective function.
Consequently, we obtain the equivalent reformulation of the original problem.
We then sequentially linearize the relaxed non-convex collision avoidance constraints while minimizing the reformulated problem with an increasing penalty term. 
The algorithm is called the penalty DCA \cite{DCA_ExacePenalty} or penalty convex-concave procedure \cite{dccp_stanford}, which aims to tighten the convexified problem of the original nonconvex one. 
This paper is organized as follows: Section II and  III presents the mathematical model of the generic trajectory planning problem for a swarm of multiple UAVs. 
Section IV and V reformulate the problem into Mixed-Integer Convex Program (MICP) and DC forms, respectively. 
The numerical results of our formulations and algorithms are shown in Section VI.
Finally, Section VII concludes the paper.

\section{State-space System Modeling of a UAV}

\begin{figure}[http]
    \centering
	\includegraphics[width=0.8\linewidth]{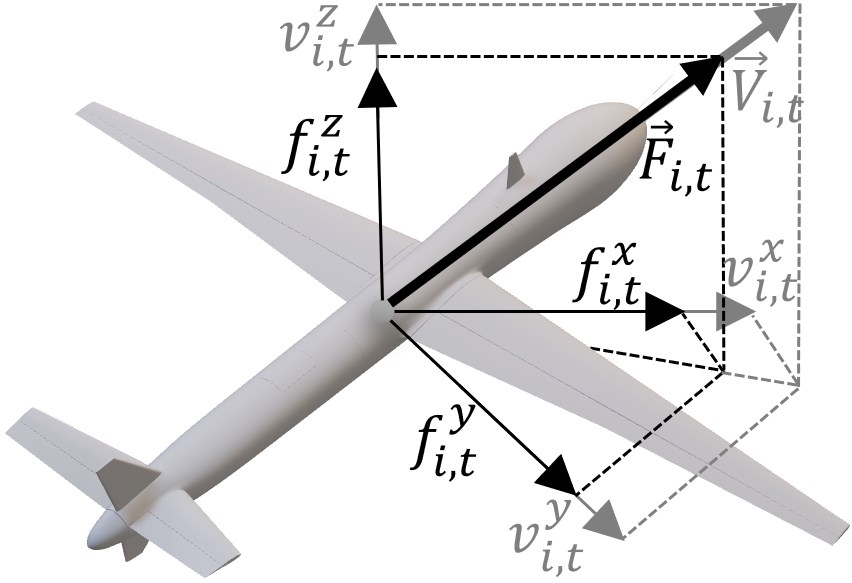}
	\vspace{-0.2pt}
	\caption{Velocity and Force vectors in body and fixed axes coordinate system}
	\label{fig: dynamics_vector}       
\end{figure}

We consider a fixed-wing UAV modeled as a point mass flying in a predetermined 3-dimensional space with the ($x,y,z$) coordinates (i.e., forward, side, and vertical directions, respectively) as shown in Figure \ref{fig: dynamics_vector} where $i$ denotes the UAV index in the swarm at the location ($x_{i,t}, y_{i,t}, z_{i,t}$) and $m_i$ is its constant mass. 
 The UAV's velocity $\Vec{V}_{i,t}$, by definition, represents the change of UAV's location as:
 \beqn
\Vec{V}_{i,t} = \Vec{v}^x_{i,t}+\Vec{v}^y_{i,t}+\Vec{v}^z_{i,t} = \frac{\mathrm{d} \Vec{x_i}}{\mathrm{d}t} 
+\frac{\mathrm{d}\Vec{y_i}}{\mathrm{d}t}
+\frac{\mathrm{d}\Vec{z_i}}{\mathrm{d}t} \label{velocity}
\eeqn
and can be decomposed into $\Vec{v}^x_{i,t} = \dfrac{\mathrm{d} \Vec{x_i}}{\mathrm{d}t}$ (the forward velocity), $v^y_{i,t}= \dfrac{\mathrm{d} \Vec{y_i}}{\mathrm{d}t}$ (the side velocity), and  $v^z_{i,t}= \dfrac{\mathrm{d} \Vec{z_i}}{\mathrm{d}t}$  (the vertical velocity). 
The force $\Vec{F}_{i,t}$ as the control input alternates the UAV acceleration following Newton's second law: 
\beqn
~ \Vec{F}_{i,t}= m_i \frac{\mathrm{d}(\Vec{V}_{i,t})}{\mathrm{d}t} = m_i\left(\frac{\mathrm{d}\Vec{v}^x_i}{\mathrm{d}t}
+\frac{\mathrm{d}\Vec{v}^y_i}{\mathrm{d}t}
+\frac{\mathrm{d}\Vec{v}^z_i}{\mathrm{d}t}\right) \label{newton}, 
\eeqn
which is also decomposed into  $\Vec{f}^x_{i,t}= \dfrac{\mathrm{d}\Vec{v}^x_i}{\mathrm{d}t}$,  $\Vec{f}^y_{i,t} = m_i \dfrac{\mathrm{d}\Vec{v}^x_i}{\mathrm{d}t}$, and $\Vec{f}^z_{i,t} = m_i \dfrac{\mathrm{d}\Vec{v}^x_i}{\mathrm{d}t}$ (i.e., forward, side, vertical forces).  

Equations (\ref{velocity})-(\ref{newton}) together form the following UAV's kyno-dynamic state-space model:
\beqn
\frac{\mathrm{d}\mathbf{x}_{i,t}}{\mathrm{d}t} = \mathbf{A} \mathbf{x}_{i,t} + B \mathbf{u}_{i,t}  \label{cont-state-space}
\eeqn
\begin{align}
\text{where}~ &
\mathbf{x}_{i,t} = \left[x_{i,t},y_{i,t},z_{i,t},v^x_{i,t},v^y_{i,t},v^z_{i,t}\right]^\top, \nonumber\\ 
&\mathbf{u}_{i,t} = \left[f^x_{i,t},f^y_{i,t},f^z_{i,t} \right]^\top, \nonumber\\
&    \mathbf{A}_i = \begin{pmatrix}
                    0  &  0  &  0   &  1  &  0   &  0\\
                    0  &  0  &  0   &  0   &  1   &  0\\
                    0  &  0  &  0   &  0  & 0  &  1 \\
                    0  &  0  &  0   &  0  &  0   & 0\\
                    0  &  0  &  0   &  0  &  0   &  0\\
                    0  &  0  &  0   &  0  &  0   &  0
                \end{pmatrix},         \nonumber
    \mathbf{B}_i = \dfrac{1}{m_i} \begin{pmatrix}
                    0   & 0      & 0\\
                      0    & 0    & 0\\
                      0    & 0                  & 0 \\
                      1 & 0      &0  \\
                      0             & 1  & 0 \\
                      0             &  0             &  1\\     
                \end{pmatrix} \nonumber 
\end{align}
Here, $\mathbf{x}_{i,t}$ denotes the vector of state variables, $\mathbf{u}_{i,t}$ denotes the control input,  $\mathbf{A}_i$ denotes the state matrix, and $\mathbf{B}_i$ denotes the input matrix of UAV $i$. 
The kyno-dynamic model (\ref{cont-state-space}) can be converted into the discrete time-variant form as follows:
   \beqn
        \mathbf{x}_{i, k+1} = \mathbf{\hat{A}}_i \mathbf{x}_{i,k} + \mathbf{\hat{B}}_i \mathbf{u}_{i,k}, ~~ \forall i \in \mathcal{N}, \forall k \in \mathcal{T} \label{statespace},
    \eeqn
where $\mathbf{\hat{A}}_i= \mathbf{I} + \Delta T \mathbf{{A}}_i,~\mathbf{\hat{B}}_i = \Delta T \mathbf{B}_i$, particularly:
    \beqn
    \mathbf{\hat{A}}_i = \begin{pmatrix}
                    1  &  0  &  0   &  \Delta T  &  0   &  0\\
                    0  &  1  &  0   &  0   &  \Delta T   &  0\\
                    0  &  0  &  1   &  0  & 0  &  \Delta T \\
                    0  &  0  &  0   &  1  &  0   & 0\\
                    0  &  0  &  0   &  0  &  1   &  0\\
                    0  &  0  &  0   &  0  &  0   &  1
                \end{pmatrix},         \nonumber
    \mathbf{\hat{B}}_i = \frac{\Delta T}{m_i}\begin{pmatrix}
                    0   & 0      & 0\\
                      0    & 0    & 0\\
                      0    & 0                  & 0 \\
                      1 & 0      &0  \\
                      0             & 1  & 0 \\
                      0             &  0             &  1\\     
                \end{pmatrix}, \nonumber 
\eeqn
and  $\mathbf{X}_{i,k} = \left[x_{i,k},y_{i,k},z_{i,k},v^x_{i,k},v^y_{i,k},v^z_{i,k}\right]^\top$ and $\mathbf{U}_{i,k} = \left[f^x_{i,k},f^y_{i,k},f^z_{i,k} \right]^\top$ respectively represent vectors of state variables and control inputs of UAV $i$ at time step $k$, and $\Delta T$ is the length of the time step.

\section{Multi-UAV Trajectory Planning Problem} 
We consider the trajectory planning problem for a swarm of $N$ UAVs in which each UAV needs to travel from its initial position to its final destination without colliding with other UAVs.  
In other words, for each UAV $i \in \mathcal{I}$ in the swarm,  we need to determine a sequence of positions ($x_{i,k}, y_{i,k}, z_{i,k}$) forming the UAV's trajectory and the sequence of control action $\mathbf{U}_{i,k}$ at each time step $k \in \mathcal{T}$ such that the UAV reaches its final destination without collision with others.
This can mathematically be formulated as a large-scale non-convex optimization problem as follows:	
\beqn
    \min_{\mathbf{x}, \mathbf{u}}  \sum_{i =1}^N \sum_{k =1}^K\bigg( \rho^f \times \underbrace{\sqrt{(f^x_{i,t})^2 + (f^y_{i,t})^2 + (f^z_{i,t})^2}}_{O_f(\mathbf{u})} \nonumber\\
  +   \rho^k \times \underbrace{\sqrt{\big( x_{i,t} - x_i^g \big)^2 + \big( y_{i,t} - y_i^g \big)^2 + \big( z_{i,t} - z_i^g \big)^2}}_{O_g(\mathbf{x})} \bigg), \label{obj_micp}
\eeqn
\noindent
subject to:
\bse
\label{individual}
\begin{align}
 &(\mathbf{x_{i}}, \mathbf{u_{i} })  \in  \Omega_{i}= \Bigg\{  \mathbf{x}_{i, k+1} = \mathbf{\hat{A}}_i \mathbf{x}_{i,k} + \mathbf{\hat{B}}_i \mathbf{u}_{i,k},~\forall k, \label{dynamics}\\
 & \quad \quad      (x_{i,1}, y_{i,1}, z_{i,1})^\top = (x_i^\mathcal{S},y_i^\mathcal{S},z_i^\mathcal{S})^\top,  \label{initial1} \\
 & \quad  \quad      \big(v^{x}_{i,1}, v^{y}_{i,1}, v^{z}_{i,1} \big)^\top = (v^{x,\mathcal{S}}_i, v^{y,\mathcal{S}}_i, v^{z,\mathcal{S}}_i)^\top,  \label{initial2}\\
        %
 & \quad \quad      (x_{i,K}, y_{i,K}, z_{i,K})^\top = (x_i^\mathcal{G},y_i^\mathcal{G},z_i^\mathcal{G})^\top, \label{goal1}\\
 & \quad  \quad      (v^x_{i,K}, v^y_{i,K}, v^z_{i,K} )^\top = (v^{x,\mathcal{G}}_i, v^{y,\mathcal{G}}_i, v^{z, \mathcal{G}}_i)^\top,
        \\
& \quad \quad         \big(f^x_{i,K}, f^y_{i,K}, f^z_{i,K} \big)^\top = (f^{x,\mathcal{G}}_i, f^{y,\mathcal{G}}_i, f^{z,\mathcal{G}}_i)^\top,
        \label{goal2}\\
&\quad \quad  \sqrt{\big(v^x_{i,k} \big)^2 + \big( v^y_{i,k} \big)^2 + \big( v^z_{i,k} \big)^2} \leq \overline{V_i},  \forall k 
        \label{velocity_conic} \\
& \quad \quad            \sqrt{\big(f^x_{i,k} \big)^2 + \big(f^y_{i,k} \big)^2 + \big( f^z_{i,k} \big)^2} \leq \overline{F_i},  \forall k 
        \label{force_conic}\Bigg\}, \forall i. 
\end{align}
\vspace{-1.5\topsep}
\ese
    \beqn
        \sqrt{( x_{i,k} - x_{j,k})^2 + ( y_{i,k} - y_{j,k})^2 + ( z_{i,k} - z_{j,k})^2} \geq d, \label{non-convex_collision} \nonumber\\
            \forall i \neq j, ~ \forall k. 
\eeqn
In the objective function \eqref{milp_collision}, we want to minimize the control effort and the traveling time of UAVs of reaching their final destination. 
 The objective consists of two terms, $O_f$ penalizes the force supplying to vehicle $i$ at time $t$ with a unit fuel cost $\rho_1$ whereas $O_g$ penalizes the remaining distance of each vehicle $i$ to its goal position multiplying with the value $\rho^k$.
 Typically, $\rho^k$ is set as an increasing function of the time indexes, e.g., $\rho^k:= a \times k,~ a>0$, so $O_g$ urges UAVs to reach their goal points as soon as possible.  The objective function is subject to two sets of constraints as follows.

 Constraint (\ref{individual}) encapsulates all local constraints state variables and control inputs of individual UAVs $i =1,\ldots, N$ in their corresponding feasible set $\Omega_i$.
In particular, the dynamics of each vehicle following discrete-time and linear state-space equation \eqref{statespace}  is now  acts as linear constraint  (\ref{dynamics}). 
The starting position is expressed  in \eqref{initial1}-\eqref{initial2} whereas the set of final conditions including the goal position, velocity, and force of vehicle are introduced in \eqref{goal1}-\eqref{goal2}.
The physical limits of UAV's velocity and driving force are captured in (\ref{velocity_conic}) and (\ref{force_conic}).
The feasible set $\Omega_i$ is convex, and (\ref{individual}) is a convex constraint. 

The constraint (\ref{non-convex_collision}) represents the collision avoidance among UAVs in the pair. 
In particular, the Euclidean separation distance between all pairs of vehicles $i \neq j$ must be equal to or greater than the safety margin $d$ at every time step $k =1, \ldots, K$.
The number of collision avoidance conditions is $\frac{N\times(N-1)}{2}\times K$.
Since the Euclidean distance norm is a convex function, (\ref{non-convex_collision}) is a non-convex constraint. 

Overall, the multi UAVs' trajectory planning problem can be summarized in the following form:
\begin{align*}
    [\mbox{P}]~~~ &\text{min}  ~ && O_f(\mathbf{u})+ O_g(\mathbf{x}) \\
        &\text{s.t.} ~ && (\mathbf{x_{i}}, \mathbf{u_{i} })  \in  \Omega_{i}, \forall i,\\
        &            ~ && \textrm{non-convex collision avoidance (\ref{non-convex_collision})}.
\end{align*}
It is worth mentioning that Problem P is generic as we can tailor $\Omega_i$ or the objective function for different application requirements, e.g., the UAV's trajectory must visit certain locations or stay close as much as possible for certain pre-determined paths.
Such modifications generally do not affect the convexity of $\Omega_i$, thus not affecting computational performance. 
The complexity of P stems from a large number of nonconvex collision avoidance conditions (\ref{non-convex_collision}).
Such constraint, however, is critical for safety requirements and cannot be ignored.   
The aim of this paper is to tackle this convex constraint set, thus facilitating the computation of UAV swarm coordination in the form of P. 

\section{Mixed-integer Convex Programming Approach}

We can use mixed integer linear programming to capture the nonconvex collision constraint (\ref{non-convex_collision}). 
Theoretically, a generic nonconvex constraint can be written in the form $x \notin \mathcal{C}$ where $\mathcal{C}$ is a convex set of variables $x$. 
If we can polyhedrally outer approximate $\mathcal{C}$ by a set of $L$ linear constraints \cite{bental_poly, glineur_SOC}
\beqn
 \mathbf{Poly}(\mathcal{C}) =\Big\{x ~ | a^\top_\nu x \leq b_\nu,~ \nu =1,2, \ldots, L \Big\}, \label{polyhedral}
\eeqn
then the condition $x \notin \mathcal{C}$ will be attained by letting at least one constraint in (\ref{polyhedral}) be violated using the auxiliary binary variable $u_m$  as follows:
\beqn
a^\top_\nu x \geq b_\nu+\epsilon - \overline{U} u_\nu,~u_\nu\in \{0,1\},~ \forall m= 1,2, \ldots M,\nonumber\\
\sum\limits_{\nu =1}^M (1-u_\nu) \geq 1 \label{micp_gen}
\eeqn
where $\overline{U}$ is a sufficient large number and $\epsilon$ is a small number. 
Constraint (\ref{micp_gen}) means that at least one value of $u_\nu=0$, consequently, one inequality $a^\top_\nu x \geq b_\nu + \epsilon$ activates, forcing $x \notin \mathcal{C}$ (the term $\epsilon$ is used to prevent the equality  $a^\top_\nu x = b_\nu$).

We are now applying (\ref{micp_gen}) to the case of the collision avoidance constraint. 
Note that we can approximate the 2-D Lorentz cone:
\[
    \mathbb{L}^2 = \Bigl\{ (\hat{x}, \hat{y}, d) \in \mathbb{R}^2 \times \mathbb{R}_{+} \Big| \sqrt{\hat{x}^2 + \hat{y}^2} \leq d \Bigl\}
\]
 by the following linear inequalities of variables $\alpha, \beta$:
\bse
\label{bental}
\beqn
       \alpha_0  \geq \hat{x},~\alpha_0  \geq -\hat{x} ,~
      \beta_0  \geq  \hat{y},~  \beta_0  \geq  -\hat{y}, \label{bental1}\\
   \Big\{     \beta_{\nu+1}  \geq -\sin \left( \frac{\pi}{2^{\nu}} \right) \alpha_{\nu} + \cos \left( \frac{\pi}{2^{\nu}} \right) \beta_{\nu},
       \label{bental2}\\
        \beta_{\nu+1}  \geq \sin \left( \frac{\pi}{2^{\nu}} \right) \alpha_{\nu} - \cos \left( \frac{\pi}{2_{\nu}} \right) \beta_{\nu},~ 
         \label{bental3}
        \\
        \alpha_{\nu+1} =\cos \left( \frac{\pi}{2^{\nu}} \right) \alpha^{\nu} + \sin \left( \frac{\pi}{2^{\nu}} \right) \beta_{\nu}, \Big\}  \label{bental4} \\
        ~\nu = 0,\dots,L-1, \nonumber \\
        \alpha_{L}  \leq d,~ \beta_{L} \leq \tan{\left(\frac{\pi}{2^L}\right)} \alpha_{L},  \label{bental5}
\eeqn
\ese
The approximation (\ref{bental}) basically forms a regular $2^L$-sided polygon with $2(L+1)$ additional variables $\alpha_\nu, \beta_\nu, \nu =1,\ldots,L$ as follows:
\beqn
 \mathbf{Poly}(\mathbb{L}^2) =\Big\{ (\hat{x},\hat{y}, \alpha, \beta) \in \mathbb{R}^{2}\times \mathbb{R}^{2(L+1)}  | (\ref{bental1})-(\ref{bental5})  \Big\}, \nonumber
\eeqn

Note also that the collision condition, i.e., the distance between two UAV is less than $d$, is in the form of 3-dimension Lorentz cone $\mathbb{L}^3$
\[
    \mathbb{L}^3 = \Bigl\{ (\hat{x}, \hat{y}, \hat{z}, d) \in \mathbb{R}^2 \times \mathbb{R}_{+} \Big| \sqrt{\hat{x}^2 + \hat{y}^2 + \hat{z}^2} \leq d \Bigl\},
\]
that can be captured by two second-order cone constraints:
\[
\hat{x}^2 + \hat{y}^2 \leq \hat{w}^2, ~ \hat{w}^2 + \hat{z}^2 \leq d^2,~ \hat{w}\geq 0,
\]
each is indeed $\mathbb{L}^2$ and can be polyhedrally approximated using (\ref{bental}). 
Consequently, we can combine  (\ref{micp_gen}) and (\ref{bental}) to construct a set of MILP constraints enforcing the distance between two UAVs outside the collision range $d$.
In particular, we need to write two sets of linear constraints (\ref{bental}) associated with the polyhedral approximation $\mathbf{Poly}(\mathbb{L}^2)$ of two 2-D Lorentz cones in the standard form (\ref{polyhedral}) and then apply the MILP reformulation trick (\ref{micp_gen}).   
Due to page limitation, we omit the presentation of the general case with arbitrary $L$. 
In the special case $L=2$, we can compact the set of constraints as follows:
\bse
\label{milp_collision}
\beqn
x_{i,k} - x_{j,k} \geq d - \overline{U} u_{i,j,k}^1,
x_{j,t} - x_{i,k} \geq d - \overline{U} u_{i,j,k}^2 \\
y_{i,k} - y_{j,k} \geq d - \overline{U} u_{i,j,k}^3,
~ y_{j,k} - y_{i,k} \geq d - \overline{U} u_{i,j,k}^4 \\
 ~ z_{i,k} - z_{j,k} \geq d - \overline{U} u_{i,j,t}^5,
~ z_{j,k} - z_{i,k} \geq d - \overline{U} u_{i,j,k}^6 \\
            \sum^6_{\nu=1} u_{i,j,k}^\nu \leq 5, \forall k, \forall i\neq j. 
  \eeqn
 which enforces the distance of two UAVs $i$ and $j$ outside the cubic outerly approximating the collision sphere of radius $d$, i.e., 
 $
     |x_{i,k} - x_{j,k}| \geq d 
  ~\text{OR} ~ |y_{i,k} - y_{j,k}| \geq d 
 ~ \text{OR}  ~ |z_{i,k} - z_{j,k}| \geq d 
  \label{abs_collision}, \nonumber\\
              \forall k, \forall i \neq j. \nonumber
$($\epsilon$ in (\ref{micp_gen}) is chosen as zero since the distance $d$ satisfies the minimum requirement of safety). 
\ese

\begin{figure}[http]
    \centering
	\includegraphics[width=1\linewidth]{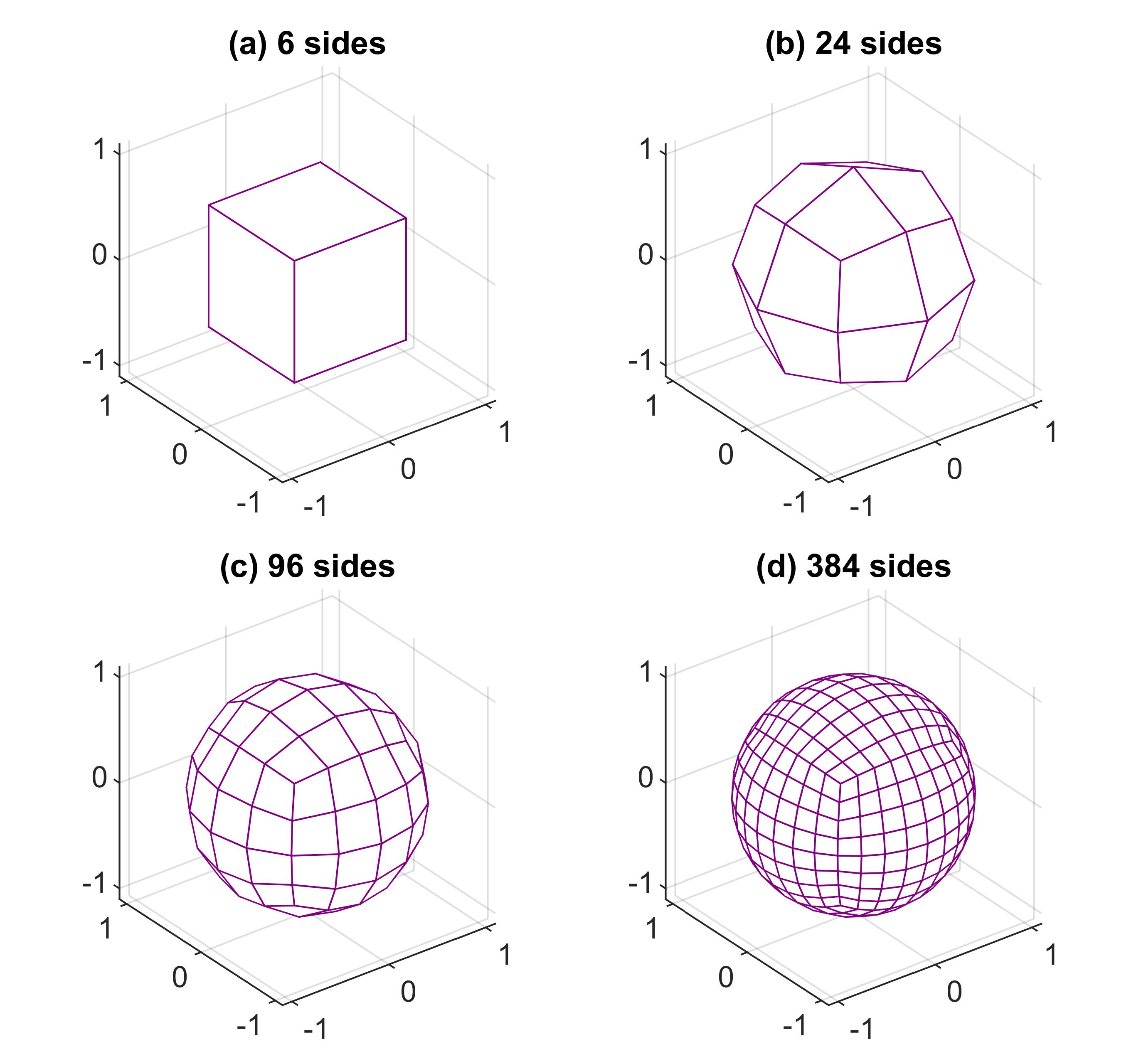}
	\vspace{-0.2pt}
	\caption{A polyhedral approximation of the 3-D ball}
	\label{fig: Cube_to_Sphere}       
\end{figure}

\noindent
\textbf{Remark:} Figure \ref{fig: Cube_to_Sphere} represents a polyhedral approximation of the 3-D ball with radius $d$.  
While the approximation error reduces as $L$ increases, the computational demand increases significantly as the number of constraints and binary variables employed increases. 
Indeed, our examination shows that only $L=2$ is computationally feasible given the number of collision avoidance conditions that we need to approximate is $\frac{N\times(N-1)}{2}\times K$.
However, MILP reformulation with $L=2$ is very conservative, which might result in infeasibility if we coordinate a large swarm of UAVs in a small space.

\section{DC Programming Approach}

\subsection{Problem Reformulation}
Let $g_{i,j,k} (\mathbf{x})$ denote the distance between two UAVs $i$ and $j$ in the time step $k$:
\beqn
    g_{i,j,k} (\mathbf{x}) =  \sqrt{( x_{i,k} - x_{j,k})^2 + ( y_{i,k} - y_{j,k})^2 + ( z_{i,k} - z_{j,k})^2}, \nonumber
\eeqn
so the collision avoidance constraints can be rewritten as 
\beqn
d - g_{i,j,k} (\mathbf{x}) \leq 0, \forall i\neq j, \forall k
    \label{dist_btw_2uav_rewrite}
\eeqn

We employ the penalty function transformation method to bring the nonconvex constraint (\ref{dist_btw_2uav_rewrite}) into an objective function of the problem P as follows:
\begin{align}
 [\text{P}_\tau]~   &\text{min}  ~ && \underbrace{O_f(\mathbf{u}) + O_g(\mathbf{x})}_{f_0(\mathbf{u,x})} + \tau\sum\limits_{i,j,k |i\neq j}  s_{i,j,k} \label{reform_obj} \\
    &\text{s.t.} ~ && (\mathbf{x}_{i}, \mathbf{u}_i)  \in  \Omega_{i}, \forall i, \label{reform_convex}\\
   &  && \left(d-g_{i,j,k}(\mathbf{x}) \right) \leq s_{i,j,k}, \forall i\neq j, \forall k \label{reform_collision}
\end{align}
where $\text{P}_\tau$ represents the penalty problem  of P with the penalty coefficient $\tau \geq 0$ and $s$ represent the relax term for original nonconvex constraint (\ref{dist_btw_2uav_rewrite}).  
There exists $\tau^* \geq 0$ such that for all $\tau \geq \tau^*$, P and $\text{P}_\tau$ have the same optimal solutions and optimal values \cite{DCA_ExacePenalty}, i.e., $s^* = 0$ and (\ref{dist_btw_2uav_rewrite}) satisfies. 
The problem $\text{P}_\tau$ is indeed a difference of the convex function (DC) programming problem, i.e., the left-hand side of (\ref{reform_collision}) can be considered as the difference of two convex functions on $\mathbf{x}$: $d$  and $g_{i,j,k}(\mathbf{x})$. 
It can be tackled by the DC Algorithm (DCA) in which we sequentially (i) solve a set of convex functions constructed by linearizing the concave term, particularly $-g_{i,j,k}(\mathbf{x})$ in (\ref{reform_collision}) (ii) increase the penalty coefficient $\tau$ until the nonconvex condition is satisfied, which will be presented next. 

\subsection{The DC Algorithm approach}
We solve the problem $\text{P}_\tau$ using the enhanced DCA, namely penalty DCA or DCA2 \cite{DCA_ExacePenalty}, or penalty convex-concave procedure \cite{dccp_stanford}, for tackling nonconvexity appearing in (\ref{reform_collision}). The algorithm is as follows:

\noindent
\textbf{Step 1:} Choose the initial point $\mathbf{\hat{x}}_0$, $\tau_0 \geq 0$, $\overline{\tau}\geq 0$, and $\mu \geq 1$.
Set the iteration $m=0$.
Initialize the set $\mathcal{O}:=\{\mathbf{\hat{x}}_0 \}$.

\noindent
\textbf{Step 2:} Solve the following optimization problem:
\begin{align}
 [\text{P}^m_\tau]~   &\text{min}  ~ && \text{objective (\ref{reform_obj})} \nonumber\\
    &\text{s.t.} ~ && \text{constraint (\ref{reform_convex})}\\
   &  && d - g_{i,j,k}(\mathbf{\hat{x}}) - \nabla^\top g_{i,j,k}(\mathbf{\hat{x}})(\mathbf{x}-\mathbf{\hat{x}}) \leq s_{i,j,k} \nonumber\\
   & &&~~~~~~~~~~~~\forall i\neq j, \forall k, \forall \mathbf{\hat{x}} \in \mathcal{O}, \label{sequential_convex} \\
   & && s_{i,j,k} \geq 0, \forall i \neq j, \forall k
\end{align}
to obtain the optimal solution $\mathbf{x}^*$.
Mathematically, we replace (\ref{reform_collision}) by a set of linear approximations at a set of points $\mathbf{\hat{x}}$ obtained so far. 

\noindent
\textbf{Step 3:} Let $\mathbf{\hat{x}}_m = \mathbf{x}^*$. Update the set $\mathcal{O}:= \mathcal{O} \cup \mathbf{\hat{x}}_m$ and update penalty coefficient $\tau_{m+1} = \min\{\mu \tau_{m}, \overline{\tau} \}$.

\noindent
\textbf{Step 4:} Stop if the following criteria satisfy:
\begin{itemize}  
\item the maximum penalty coefficient reaches
\[
\tau_m = \overline{\tau}
\]
    \item the gap between optimal objectives found between two consecutive iterations is small  
\beqn
    \begin{aligned}
        \delta_m = \Bigg( f_0(\mathbf{u}^*_m, \mathbf{x}^*_m)  + \tau_m \sum\limits_{i,j,k|i\neq j} s^*_{i,j,k,m} \Bigg) - \nonumber \\
        \Bigg( f_0(\mathbf{u}^*_{m-1}, \mathbf{x}^*_{m-1})  + \tau_{m-1} \sum\limits_{i,j,k|i\neq j} s^*_{i,j,k,m-1} \Bigg) \leq \epsilon , \nonumber
    \end{aligned}
\eeqn
\end{itemize}
where $\epsilon$ is a very small number acting as  the tolerance.
Note also that $\mathbf{u}^*_m, \mathbf{x}^*_m, s^*_{i,j,k,m}$ are optimal solutions of $\mathbf{u}, \mathbf{x}, s$ found by solving $\text{P}^m_\tau$. If the stopping conditions are not satisfied, update $m= m+1$ and go back to \textbf{Step 1}.

The iterative algorithm consists of 4 steps.  
The key point is that for each iteration we replace the distance between UAV $i$ and $j$ at time $k$ by its linearization at $\mathbf{\hat{x}}_m$, 
\[
g_{i,j,k}(\mathbf{x}) := g_{i,j,k}(\mathbf{\hat{x}}_m) + \nabla^\top g_{i,j,k}(\mathbf{\hat{x}}_m)(\mathbf{x}-\mathbf{\hat{x}}_m) 
\]
and consequently obtain the linear approximation of (\ref{reform_collision}) at $\mathbf{\hat{x}}_m$ as follows:
\[
d - g_{i,j,k}(\mathbf{\hat{x}}_m) - \nabla^\top g_{i,j,k}(\mathbf{\hat{x}}_m)(\mathbf{x}-\mathbf{\hat{x}}_m) \leq s_{i,j,k} \forall i\neq j, \forall k
\]
Consequently, we obtain the convex optimization problem $\text{P}^m_\tau$ in Step 2.
Over iterations, the set of linearized constraints (\ref{sequential_convex}) expands to tighten the convexification of the constraint (\ref{reform_collision})  whereas the increasing $\tau^m$ due to $\mu >1$ enforce the slack variables $s$ converge to zero.
Together, they try to enforce the feasibility of the obtained solution, i.e., the nonconvex collision avoidance (\ref{dist_btw_2uav_rewrite}) satisfy and the optimal values of $\text{P}^m_\tau$ converge to the sub-optimal values of P.
In other words, we aim to obtain an upper bound of P with a feasible solution $\mathbf{x}^*$.

\noindent
\textbf{Remark:} Unlike the MICP formulation, which is NP-hard, the DC programming approach enables us to solve the UAV planning problem by sequentially solving a set of convex program $\text{P}_\tau^m$.
As each convex program can be solved efficiently by matured convex optimization algorithms such as interior point methods, the computational performance can be improved significantly.
Mathematically, MICP requires approximating the non-convex feasible set (\ref{non-convex_collision}) beforehand by employing a set of a large number of MILP constraints (\ref{milp_collision}.
Many constraints in this set are non-binding at optimum and can be ignored. 
In contrast, in DC programming, we sequentially add the linearization of the nonconvex constraints at explored points found after each iteration.

\section{Numerical Results}


We implemented DC programming approach on a PC configured with an Intel Xeon and 32GB of RAM. 
To benchmark the performance of both models, we verify their formulation for 5, 10, and 15 vehicles with the GUROBI solver.
Consequently, the number of collision avoidance conditions needed to be satisfied at each time step is 10, 45, and 105.
In the three numerical experiments, the minimum safety distance between vehicles is $d=5$ distance units, and the quantity of time steps is $\mathcal{T} = 30$ time units. 
We compare the DC programming results with the ones obtained by using MICP model with the cubic approximation (\ref{abs_collision}) of collision avoidance. 

\begin{figure}[http]
    \centering
	\includegraphics[width=0.8\linewidth]{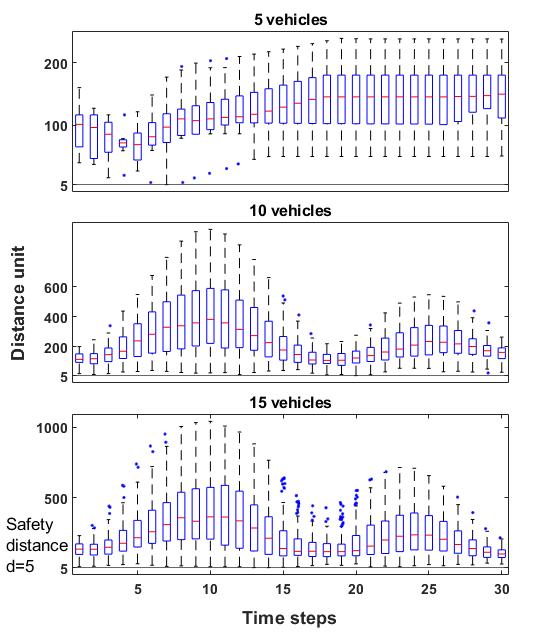}
	\vspace{-15pt}
	\caption{Distance between vehicles DCA Model}
	\label{fig: DCA_Model}       
\end{figure}

Fig. \ref{fig: DCA_Model} shows results of the distance between vehicles at each time step $k$ obtained by solving the UAV planning problem using DC programming approaches.
It shows that there is no crash between vehicles throughout the time steps in the DC model in all three experiments.
In other words, the DC programming approach guarantees the satisfaction of a large number of  nonconvex collision avoidance conditions.   

\begin{figure}[http]
    \centering
	\includegraphics[width=0.75\linewidth]{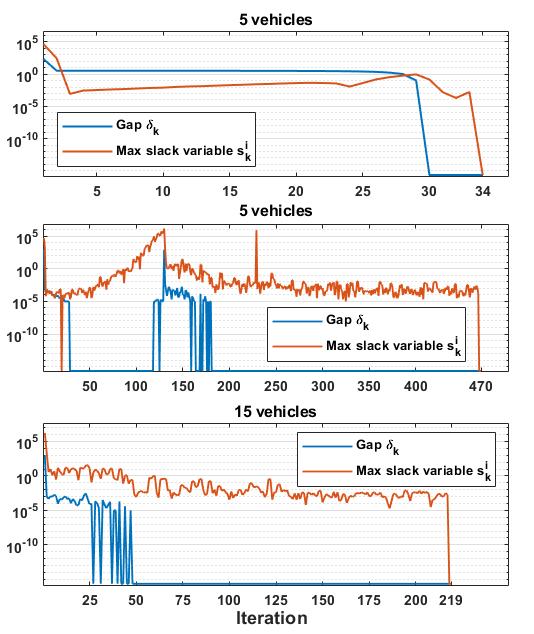}
	\vspace{-10pt}
	\caption{Convergence analysis of DCA for experiments}
	\label{fig: Covergence_15UAV}       
\end{figure}

Fig. \ref{fig: Covergence_15UAV} demonstrates the numerical convergence for penalty DCA used to solve the DC programs in all test cases.
The maximum value among all slack variables $s^*_{i,j,k}\geq 0$ converge to zero, which means all collision avoidance constraints are also satisfied at the optimum and also the objective value is equal to the original one, i.e., the penalty term $\tau \sum\limits_{i,j,k|i\neq j} s^*_{i,j,k} =0$.
Additionally, the gap between the objective function found between two consecutive iterations $\delta_m$ converges to zero, which means we reach the local optimum (sub-optimal solution) is found. 
In our experiment, the optimal solutions of 5-vehicle, 10-vehicle, and 15-vehicle experiments are converged at iterations 34, 470, and 219, respectively.

\begin{figure}[t!]
\centering
\subfigure[Three-dimensional view]{\includegraphics[width=0.6\linewidth]{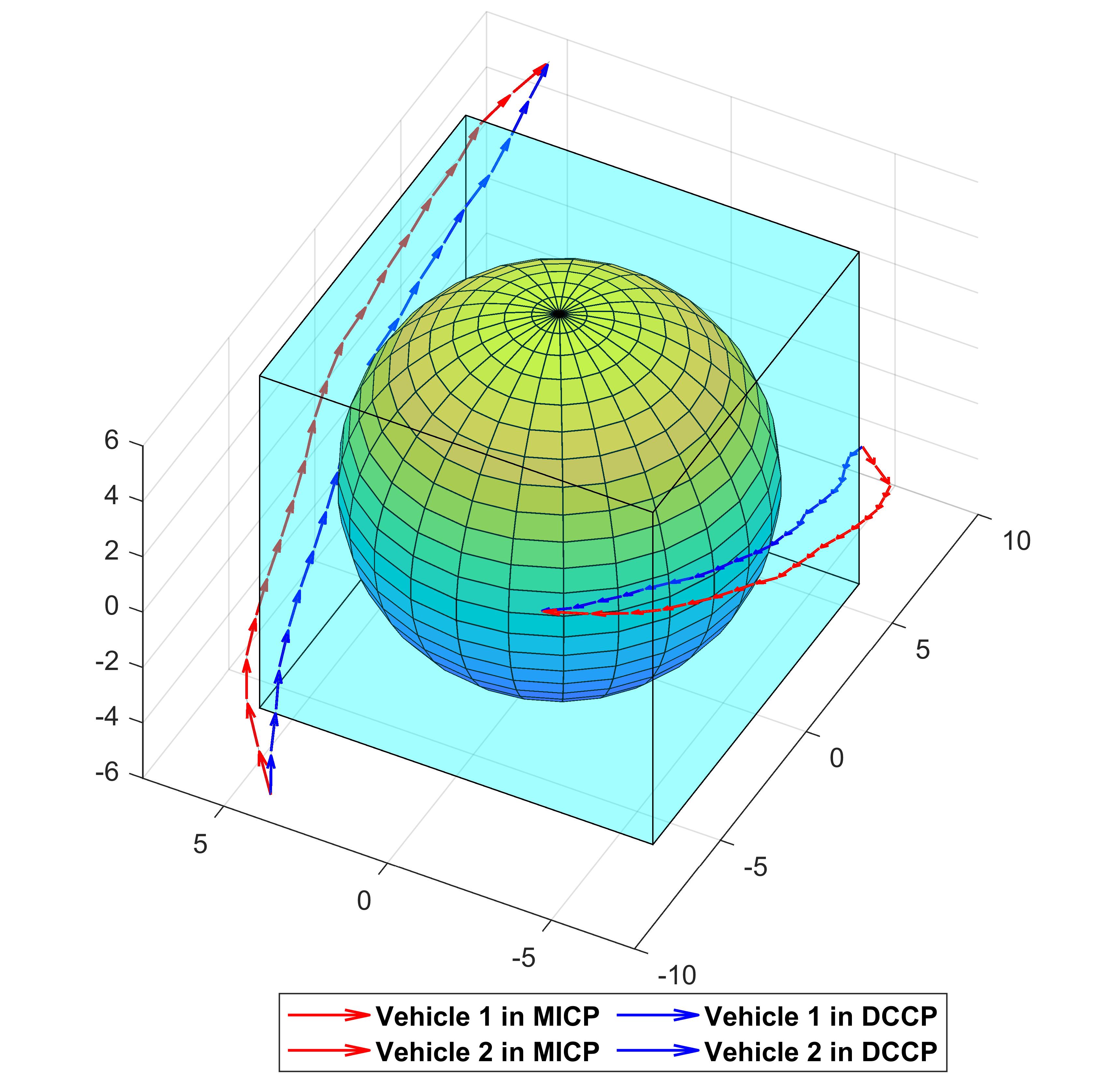}} \\
\vspace{-1.5\topsep}
\subfigure[Top view]{\includegraphics[width=0.6\linewidth]{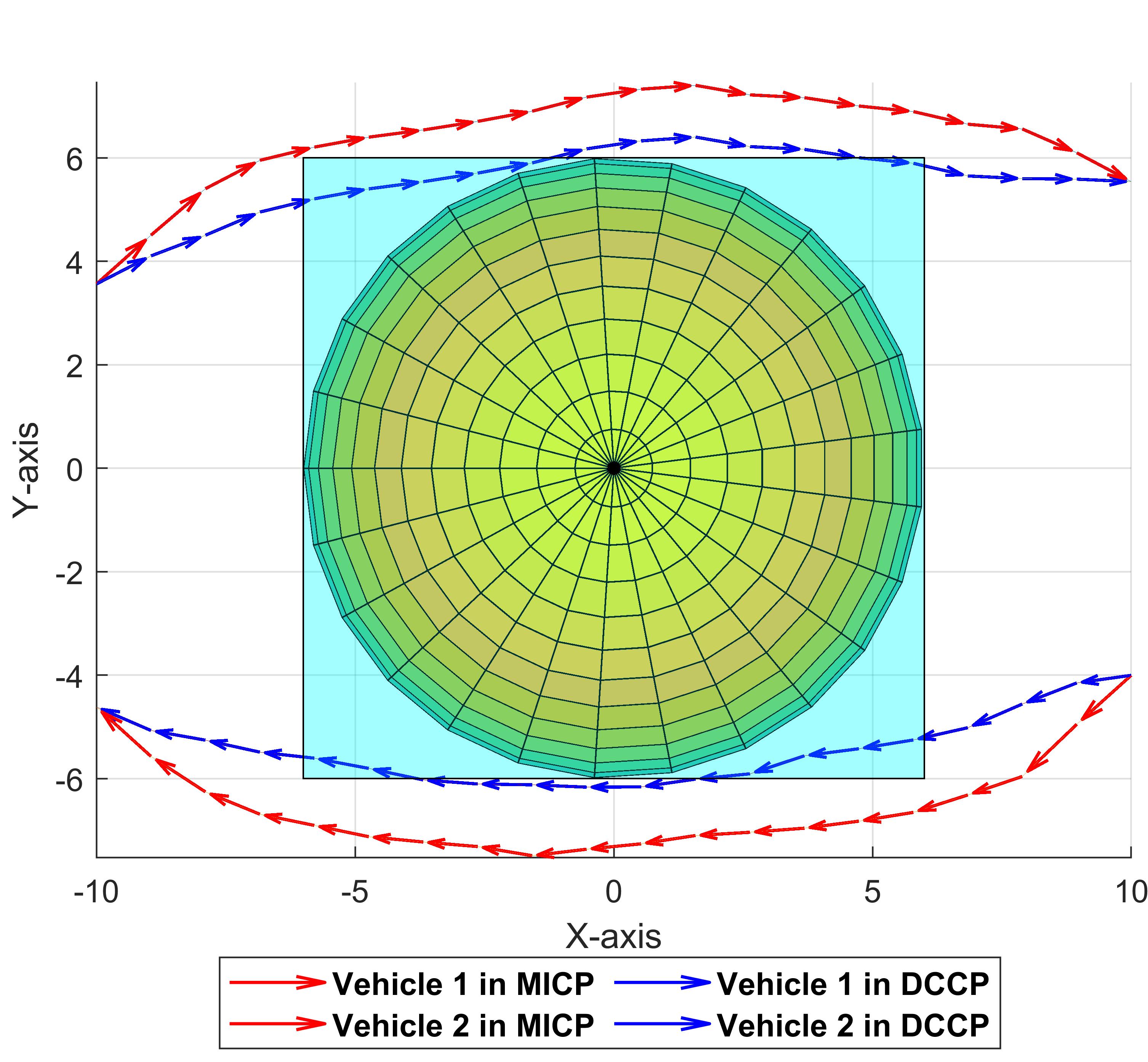}}%
\vspace{-5pt}
\caption{Illustration of the collision avoidance constraint in MICP model and DCA model for trajectory planning of two vehicles}
\label{fig: collision_avoidance}
\end{figure}

The obtained sub-optimal solution of DC program generally has a very good performance, even surpassing the MICP approach. 
This is because the DC programming approach employs a less conservative approximation of the nonconvex collision condition, as shown in Fig. \ref{fig: collision_avoidance}. 
In the DCA model, we can utilize the full collision-free space outside the radius sphere $d$ (safety distance).
In contrast, the  cubic approximation (\ref{abs_collision}) used in the MICP is more conservative.
Therefore, the fuel cost of DC model is lower than that of MICP, as shown in the Table. \ref{tab:Fuel_cost}.
Note also that, while increasing the size of the polyhedral approximation (as shown in Fig. \ref{fig: Cube_to_Sphere}) can reduce the conservatives, the MICP easily becomes intractable. 
Indeed, only the cubic approximation (\ref{abs_collision}) \cite{taskallocation2001} widely used in the literature is computationally feasible in our experiments.

\vspace{-10pt}

\begin{table}[http]
\centering
\footnotesize
\caption{Grand total fuel cost of vehicles in MICP and DCA}
\vspace{-0.2cm}
\label{tab:Fuel_cost}
\begin{tabular}{lccc}
\hline
 & \begin{tabular}[c]{@{}c@{}} \textbf{MICP} \end{tabular} & \begin{tabular}[c]{@{}c@{}} \textbf{DCP} \end{tabular} & \begin{tabular}[c]{@{}c@{}} \textbf{$\Delta ($ MICP $-$ DCP $)$} \end{tabular} \\ \hline
5 vehicles & 326.18   & 323.17  & 3.01 \\ \hline
10 vehicles & 1594.11 & 1592.32 & 1.79 \\ \hline
15 vehicles & 2855.25 & 2839.97 & 15.28 \\ \hline
\end{tabular}
\end{table}

\section{Conclusion}
This paper examines the use of the DC programming approach to solve the planning problem of a UAV swarm considering the nonconvex collision avoidance requirement. 
In particular, we sequentially approximate this nonconvex constraint by its linearization and adopt the penalty reformulation with slack variables. 
The problem is effectively tackled by sequentially solving a set of computationally manageable convex programs. 
Compared to the traditional mixed integer optimization model with the cubic approximation of collision avoidance constraint, the obtained solution satisfies the safety condition while achieving better cost saving thanks to its less conservative approach. 





\bibliographystyle{IEEEtran}
\bibliography{ref}

\end{document}